\documentclass[11pt]{article}

\usepackage{arxiv}

\usepackage[utf8]{inputenc} % allow utf-8 input
\usepackage[T1]{fontenc}    % use 8-bit T1 fonts
\usepackage{hyperref}       % hyperlinks
\usepackage{url}            % simple URL typesetting
\usepackage{booktabs}       % professional-quality tables
\usepackage{amsfonts}       % blackboard math symbols
\usepackage{nicefrac}       % compact symbols for 1/2, etc.
\usepackage{microtype}      % microtypography
\usepackage{lipsum}
\usepackage{hyperref}
\usepackage[all]{xy}
\usepackage{amssymb,amsmath}

\def\sympow{{\setbox0\hbox{$\bigcirc$}\setbox1\hbox to\wd0{\hss$s$\hss}%
\wd1 0pt\box1\box0}}%symmetric power

\newcommand{\pproof}[1]{\textbf{Proof. }{#1} $\blacksquare$}
\newcommand{\Nset}{\mathbb{N}}

\newcommand{\Rset}{\mathbb{R}}
\newcommand{\Zset}{\mathbb{Z}}
\newcommand{\Cset}{\mathbb{C}}

\newcommand{\Identidad}{{\rm Id}}

\newcommand{\rme}{{\rm e}}
\newcommand{\rmi}{{\rm i}}
\newcommand{\Orden}{\mathop{\rm O}\nolimits}
\newcommand{\ds}{\displaystyle}

\def\WKB{K_{\textsf{WKB}}}
\usepackage[colorinlistoftodos]{todonotes}
\usepackage[normalem]{ulem}

\definecolor{ColorTomasL}{rgb}{0.9,0.2,0.1}

\definecolor{ColorJJ}{rgb}{0.9,0.,0.4}

\definecolor{ColorPrimi}{rgb}{0.1,0.5,0.1}

\definecolor{ColorChara}{rgb}{0.1,0.2,0.9}

\newcommand{\A}{\mathcal{A}}

\newcommand{\xbf}{\mathbf{x}}
\newcommand{\erfi}{\mathrm{erfi}}
\newcommand{\erf}{\mathrm{erf}}
\newcommand{\hot}{\mathrm{h.o.t.}}

\newtheorem{corollary}{Corollary}

\newtheorem{definition}{Definition}

\newtheorem{proposition}{Proposition}

\usepackage{verbatim}

\usepackage{cite}
\usepackage{mathrsfs}
\title{Semiclassical quantification of some two degree of freedom potentials: a Differential Galois approach}

\author{Primitivo Acosta-Hum\'anez \\
Instituto \& Escuela de Matem\'atica \\
Universidad Aut\'onoma de Santo Domingo\\
Dominican Republic \\
  \texttt{pacosta-humanez@uasd.edu.do} \\
   \And
 J. Tom\'as L\'azaro\\
 Departament de Matem\`atiques\\
 Universitat Polit\`ecnica de Catalunya (UPC)\\
 Centre de Recerca Matem\`atica (CRM) \\
 Institute of Mathematics of the UPC-BarcelonaTech (IMTech)\\ Barcelona, Spain\\
  \texttt{jose.tomas.lazaro@upc.edu} \\
   \AND
  Juan J. Morales-Ruiz\\
 Departamento de Matem\'atica Aplicada\\
 Universidad Polit\'ecnica de Madrid (UPM)\\
 Madrid, Spain\\
  \texttt{juan.morales-ruiz@upc.edu} \\
   \And
   Chara Pantazi\\
 Departament de Matem\`atiques\\
 Universitat Polit\`ecnica de Catalunya (UPC)\\
 Barcelona, Spain\\
  \texttt{chara.pantazi@upc.edu} \\
}
\begin{document}
\maketitle
\begin{abstract}
In this work we explain the relevance of the Differential Galois Theory in the semiclassical (or WKB) quantification of some two degree of freedom potentials. The key point is that the semiclassical path integral quantification around a particular solution depends on the variational equation around that solution: a very well-known object in dynamical systems and variational calculus. Then, as the variational equation is a  linear ordinary differential system,  it is possible to apply the Differential Galois Theory to study its
solvability in closed form. We obtain closed form solutions for the semiclassical quantum fluctuations around constant velocity solutions for some systems like the classical Hermite/Verhulst, Bessel, Legendre, and Lam\'e potentials.
%Verhulst’s potentials.
We remark that some of the systems studied are not integrable, in the Liouville - Arnold sense.
\end{abstract}

% keywords can be removed

\section*{Introduction}
In~\cite{m2020} the third author suggested the relevance that the Differential Galois Theory could play in the Feynman's path integral approach in Quantum Mechanics. Indeed, the key proposal was to study whether it was possible to obtain, in closed form, the semiclassical approximation of the Feynman's propagator $K$, see \cite{pauli2000,cecile1976}.

Let us recall the notation and main ideas in~\cite{m2020}.  Given $n$ the number of degrees of freedom, we denote by
${\bf x}=(x_1,x_2,\ldots,x_n)$ the position, %$y=(y_1,y_2,\ldots,y_n)=\dot{x}$ momenta,
$t$ the time, and $\gamma$ is a path from  $({\bf x}_0,t_0)$ to $({\bf x}_1,t_1)$. This classical path $\gamma$ in the configuration space defines an integral curve $\Gamma$ in the phase space, assuming there are non focal (conjugated) points (chapter 9 of \cite{arnol2013}).

The computation of the  propagator $K({\bf x}_1,t_1 \,|\, {\bf x}_0,t_0)$ around the path $\gamma$ in the semiclassical approach (where $\hslash$ is small) can be obtained through
\[
K({\bf x}_1,t_1 \,|\, {\bf x}_0, t_0) = K_{\textsf{WKB}} \left(1 + \Orden(\hslash)\right),
\]
 ${K_\textsf{WKB}}$ being the semiclassical approximation of the propagator $K$ (WKB after Wentzel, Kramers, and Brillouin, 1926). The function  $\WKB$ is given by the so-called \emph{Pauli-Morette formula}
\begin{equation}\label{Pauli}
K_{\textsf{WKB}}({\bf x}_1,t_1 \,|\, {\bf x}_0,t_0) =A\, \rme^{\frac{\rmi}{\hslash} S(\gamma)},\qquad \textrm{with} \qquad  A=\frac{1}{\left( 2\pi i\hslash\right)^{n/2}} \frac{1}{\sqrt{\det J(t_1,t_0)}},
\end{equation}
\ \
called prefactor, where:
\begin{itemize}

\item The {\it fixed classical path $\gamma$} parametrized by $({\bf x}(t),t)$, for $t_0\leq t \leq t_1$, with starting point $({\bf x}_0,t_0)$ and endpoint $({\bf x}_1,t_1)$, assuming no focal points.

\item $S$ is the \emph{action} computed on this classical path $\gamma$, i.e.
\begin{equation}
S[{\bf x}(t)]:=S(\gamma) = \int_{t_0}^{t_1} L({\bf x},\dot{{\bf x}},t) \, dt =
\int_{t_0}^{t_1} \left( \sum_{i=1}^n y_i \dot{x}_i - H({\bf x},{\bf y},t) \right) \, dt,
\label{def:action}
\end{equation}
with $L({\bf x},\dot{{\bf x}},t)$ and $H({\bf x},{\bf y},t)$ being the corresponding \emph{Lagrangian} and \emph {Hamiltonian} functions.

\item The
$n\times n$ matrix  $J=J(t_1,t_0)$  is given by
a block inside a fundamental matrix $\Psi (t,t_0)$ of the variational equations around the phase integral curve defined by the classical path $\gamma$:
\[
\Psi (t,t_0)=\begin{pmatrix} \quad \circ &J(t,t_0)\\ \quad \circ & \circ \end{pmatrix}.
\]
That is, if the fundamental matrix (with initial condition $\Psi(t_0,t_0)= \mathrm{Id}_{2n}$, the $2n$-dimensional identity matrix) is splitted into four square boxes of dimension $n$, the matrix $J$ is given by the variation of the positions with respect to the initial momenta.
We will refer to $J$ and $\det J$ as the \emph{Van Vleck-Morette matrix and its determinant}, respectively.
\end{itemize}

We consider  the above  connection between the semiclassical propagator and the variational equation as a quantum mechanical confirmation  of the following fundamental {\it Bryce DeWitt's Principle}~\cite{DeWitt1965}:
\begin{center}
\emph{``The quantum theory is basically a theory of small disturbances"}
\end{center}
So, the computation of the semiclassical propagator in formula~\eqref{Pauli} is based on the matrix $J$, obtained in its turn  by means of the solution of the variational equation. Recall that the variational equation is also called Jacobi equation (in the context of variational calculus), equation of geodesic dispersion (in general relativity) or equation of small disturbances (in agreement with Bryce DeWitt's Principle).  It is easy to see that
\begin{equation}
J(t_1,t_0)= \left( \frac{\partial \xi_i(t_1)}{\partial \eta_j(t_0)} \right),
\label{def:J}
\end{equation}
where $(\xi_1....,\xi_n, \eta_1,...,\eta_n)$ are the variables in the corresponding variational equation, and $\xi_i=\delta x_i$, $\eta_i=\delta y_i$ being the variations in positions and momenta, respectively.  Since the variational equation is a linear differential system it is possible to study its solutions by means of the differential Galois theory (see~\cite{morales2013} for details).
In the case that the classical Hamiltonian system under consideration is integrable in \emph{Liouville} - \emph {Arnold} sense then one of the main results in~\cite{m2020} guarantees that it is possible to obtain a closed form formula for $\WKB$, in a very precise way. Keep in mind, however,
that this is only a necessary condition for integrability.
Indeed, as it will be seen, some of the Hamiltonian families considered along this paper are not integrable {\it but admit closed analytic formulas for the semiclassical propagator} around some special classical paths.

Along this paper integrability of the variational equation means integrability in the sense of the differential Galois theory. This kind of integrability is characterized by the structure of the  Galois group of the equation: the identity component of the Galois group must be solvable (see the appendix of~\cite{m2020} and references therein).

One of the open problems posed in~\cite{m2020} was to compute this propagator for concrete families of Hamiltonian systems with several degrees of freedom. This paper aims to be a first work filling this gap.

To do that, we consider $2$-degrees of freedom Hamiltonian systems and compute the propagator along classical paths $\gamma$ defined by invariant planes in the phase space, whose restricted dynamics fall in a free particle model of $1$-degree of freedom. That is, defined by
\[
E=\frac{p_x^2}{2},\qquad \qquad p_x=\dot x,
\]
whose solution curve is given by $(x_E(t),p_{x_E}(t))$, where
\[
x_E(t)=\sqrt{2E}t+x_0, \qquad \qquad p_{x_E}(t)=\sqrt{2E}.
\]
Since only autonomous Hamiltonian systems are considered we can, without loss of generality, take initial time $t_0=0$. Therefore, for $n=2$ the Pauli-Morette formula~\eqref{Pauli} reads
\begin{equation}\label{Paulia}
\WKB(\xbf_1,t_1 \,|\,  \xbf_0,0) = \frac{1}{ 2\pi \rmi\hslash}
\frac{1}{\sqrt{\det J(t_1,0)}} \, \rme^{\frac{\rmi}{\hslash} S(\gamma)}.
\end{equation}
We want to stress that we shall only semiclassically quantify around this very special type of solutions, i.e. those of a free particle motion, their most important quantum oscillations occurring in their transversal directions. We do not know the physical relevance of this kind of quantification around such a path.
In any case, to quantify around particular special solutions is still today a very frequent method in quantum mechanics and in quantum field theory. For instance, the case around \emph{instantons} in tunneling problems.

Furthermore, is very likely that  for the families of potentials considered here it would be impossible to quantify by means of closed form formulas around an arbitrary classical path. The reasons are clear: before quantifying,  it is necessary to get an analytic expression for the particular integral curve of the classical mechanical system (or the associate classical path).
%we have to obtain analytically the  particular integral curve of the classical mechanical system (or the associate classical path) and it is not clear at all how to do that for the families of potentials considered here.
Indeed, as mentioned above, some of the families considered here are not integrable,  and hence, such general solution in closed form does not exist. Thus, we state the following results (for further details see Section \ref{finrem}):

\textit{Consider the families of potentials defined by equation \eqref{def:potential} with $k=2$ given by Table \eqref{table}. Then, the Bessel and Legendre families with $b\neq 0$ of Hamiltonian systems, as well as Hermite and Lam\'e families are not integrable in the Liouville-Arnold sense. Furthermore the semiclassical approximations of the corresponding Feynman propagators are not integrable.}

We remark that closed form solutions, also in the framework of the differential Galois theory, of some  non-autonomous one-dimensional  oscillators were obtained by the first named author in \cite{as2013, akss2015,as2015}.

The paper is structured as follows: Section \ref{se:2dofhamiltonians} is devoted to some generalities regarding the structure of the Van Vleck-Morette determinant for $2$-degrees of freedom Hamiltonian systems around invariant planes with free particle reduced dynamics. In Section~\ref{appl} we obtain closed analytical formulas for the Van Vleck-Morette determinant and hence the semiclassical approximation of the Feynman propagators for four classical families of potentials:
Hermite, Bessel, Legendre and Lam\'e, see Table \ref{table}. For any family, an illustrative example accompanying the theoretical result is provided. Finally, in Section \ref{finrem} we prove Proposition \ref{prop:noint} and Corollary \ref{GOAL}.

%%%%%%%%%%%%%%%%%%%%%%%%%%%%%%%%%%%%%%%%%%%%%%%%%%%%%%%%%%

\section{Potentials with free particle motion in an invariant plane}
\label{se:2dofhamiltonians}

Let us consider classical Hamiltonian systems with $n=2$ degrees of freedom. Assume that they are given by the sum of a kinetic energy $T$ and a certain kind of potentials $V$. More precisely,
if $(x,y)$ stands for the spatial variables and $(p_x,p_y)$ for the corresponding momenta, then the Hamiltonian function is $H=H(x,p_x,y,p_y)=T(p_x,p_y)+V(x,y)$, where
\[
T(p_x,p_y) =\frac{p_x^2+p_y^2}{2},  %\frac{\dot{x}^2+\dot{y}^2}{2},
\]
and
\begin{equation}
V(x,y)=y^k f(x,y), \qquad \quad f(x,0)\ne 0, \quad k\in \Nset, \ k\geq 2.
\label{def:potential}
\end{equation}
Despite of its concrete form, this kind of potentials has been considered in references as~\cite{morales2013,ab2008} among others.
The system of ordinary differential equations associated to $H(x,p_x,y,p_y)$ reads
\begin{equation}
\begin{array}{lclclcl}
\dot{x} &=&\dfrac{\partial H}{\partial p_x} = p_x  & \qquad \qquad &
\dot{y} &=&\dfrac{\partial H}{\partial p_y} = p_y \\[1.9ex]
\dot{p}_x &=& -\dfrac{\partial H}{\partial x} = -y^k \dfrac{\partial f}{\partial x}(x,y) & \qquad \qquad &
\dot{p}_y &=& -\dfrac{\partial H}{\partial y} = -ky^{k-1} f(x,y) - y^k \dfrac{\partial f}{\partial y}(x,y) \\[1.9ex]
\end{array}
\label{def:general:system}
\end{equation}
where dot means derivative with respect to $t$. This kind of Hamiltonian systems has, for $k\geq 2$, the invariant plane $\Gamma=\{y=p_y=0\}$ (in the phase space), with associated
$\{y=0\}$ invariant straight line path $\gamma$ in the configuration space.

We restrict ourselves to solutions of~\eqref{def:general:system} lying on $\Gamma$, which connect two points $(x_0,p_{{x}_0},0,0)$ and $(x_1,p_{{x}_1},0,0)$, at times $t_0=0$ and $t=t_1$, respectively,
and such that its motion is the one of a free particle, that is,
\begin{equation}
x(t)=\dfrac{x_1-x_0}{t_1}\, t+x_0, \qquad \qquad p_x(t)= \dfrac {x_1-x_0}{t_1}.
\label{param:free:particle}
\end{equation}
Fixed an energy level $H(x,p_x,y,p_y)=E$, we denote by $\gamma_E$ the free particle path~\eqref{param:free:particle} associated to this energy $E$. This implies, in particular, the relation
$E=\frac{1}{2}p_x^2$.
%$p_x=\sqrt{2E}%
Hence,
\[
E = H(x,\dot{x},0,0)= \frac{1}{2}\dot{x}^2 \Longrightarrow \dot{x}(t)= \sqrt{2E} \Longrightarrow x(t) = \sqrt{2E} \, t + x(0),
\]
it follows that the (free-particle) path $\gamma_E$ can be parameterized by

\begin{equation}
x_E(t):=\sqrt{2E} \, t + x_0, \qquad x_0=x(0).
\label{freeparticle:xH}
\end{equation}
Moreover,
\begin{equation*}
E=\frac{1}{2} \left( \frac{x_1-x_0}{t_1} \right)^2 >0,
\label{energy:x0:x1}
\end{equation*}
from which the action on $\gamma_E$ becomes
\begin{eqnarray*}
S[\gamma_E] &=& \int_{t_0=0}^{t_1} \left( p_x \dot{x} + p_y \dot{y} - H(x,p_x,y,p_y) \right) \, dt
= \int_{0}^{t_1} \left( p_x^2(t) - \frac{1}{2}p_x^2(t) \right) \, dt \\
&=&
\int_{0}^{t_1} E \, dt = E t_1 = \frac{1}{2} \frac{(x_1-x_0)^2}{t_1}.
\end{eqnarray*}
To avoid misunderstandings, henceforward we will fix the following order $(x,p_x,y,p_y)$ for the variables, and will denote by $X_H=(F_1,F_2,F_3,F_4)$ the right hand-side vector field in~\eqref{def:general:system}.
In general, if $\Gamma$ is any solution of~\eqref{def:general:system} then the corresponding \emph{variational equation} around $\Gamma$ is defined as
\[
\left(
\begin{array}{c}
\dot{\xi}_1 \\[1.2ex]
\dot{\eta}_1 \\[1.2ex]
\dot{\xi}_2 \\[1.2ex]
\dot{\eta}_2
\end{array}
\right) = \left(
\begin{array}{cccc}
\dfrac{\partial F_1}{\partial x} & \dfrac{\partial F_1}{\partial p_x} &
\dfrac{\partial F_1}{\partial y} & \dfrac{\partial F_1}{\partial p_y} \\[1.9ex]
\dfrac{\partial F_2}{\partial x} & \dfrac{\partial F_2}{\partial p_x} &
\dfrac{\partial F_2}{\partial y} & \dfrac{\partial F_2}{\partial p_y} \\[1.9ex]
\dfrac{\partial F_3}{\partial x} & \dfrac{\partial F_3}{\partial p_x} &
\dfrac{\partial F_3}{\partial y} & \dfrac{\partial F_3}{\partial p_y} \\[1.9ex]
\dfrac{\partial F_4}{\partial x} & \dfrac{\partial F_4}{\partial p_x} &
\dfrac{\partial F_4}{\partial y} & \dfrac{\partial F_4}{\partial p_y} \\[1.9ex]
\end{array}
\right)_{\bigg|_{\Gamma}}
\left(
\begin{array}{c}
\xi_1 \\[1.2ex]
\eta_1 \\[1.2ex]
\xi_2 \\[1.2ex]
\eta_2
\end{array}
\right),
\]
where $\xi_j$ stands for the positions and $\eta_j$ for the momenta.
In the case of system~\eqref{def:general:system}, the variational equation around any solution lying on $\Gamma$ takes the form
\begin{equation}
\left(
\begin{array}{c}
\dot{\xi}_1 \\
\dot{\eta}_1 \\
\dot{\xi}_2 \\
\dot{\eta}_2
\end{array}
\right) = \left(
\begin{array}{cccc}
0 & 1 & 0 & 0 \\
\A_1(x,y) & 0 & \A_2(x,y) & 0  \\
0 & 0 & 0 & 1 \\
\A_2(x,y) & 0 &
\A_3(x,y) & 0
\end{array}
\right)
\left(
\begin{array}{c}
\xi_1 \\
\eta_1 \\
\xi_2 \\
\eta_2
\end{array}
\right),
\label{VE:y:py:zero:general-case}
\end{equation}
with
\begin{eqnarray*}
\A_1(x,y)&=&-y^k \frac{\partial^2 f}{\partial x^2}(x,y), \qquad
\A_2(x,y)= -k y^{k-1} \frac{\partial f}{\partial x}(x,y) - y^k \frac{\partial^2 f}{\partial x \partial y}(x,y), \\[1.4ex]
\A_3(x,y)&=&-k(k-1)y^{k-2}f(x,y) - 2ky^{k-1} \frac{\partial f}{\partial y}(x,y) - y^k\frac{\partial^2 f}{\partial y^2}(x,y),
\end{eqnarray*}
for $k\geq 2$.

On one hand, the case $k>2$ becomes
\begin{equation*}
\left(
\begin{array}{c}
\dot{\xi}_1 \\
\dot{\eta}_1 \\
\dot{\xi}_2 \\
\dot{\eta}_2
\end{array}
\right) = \left(
\begin{array}{cccc}
0 & 1 & 0 & 0 \\
0 & 0 & 0 & 0 \\
0 & 0 & 0 & 1 \\
0 & 0 & 0 & 0
\end{array}
\right)
\left(
\begin{array}{c}
\xi_1 \\
\eta_1 \\
\xi_2 \\
\eta_2
\end{array}
\right),
\end{equation*}
with straightforward solutions.

On the other, the case $k=2$ is much richer and is, therefore, the one tackled in this work.

Indeed, for $k=2$ variational equation reduces to
\begin{equation}
\left(
\begin{array}{c}
\dot{\xi}_1 \\
\dot{\eta}_1 \\
\dot{\xi}_2 \\
\dot{\eta}_2
\end{array}
\right) = \left(
\begin{array}{cccc}
0 & 1 & 0 & 0 \\
0 & 0 & 0 & 0 \\
0 & 0 & 0 & 1 \\
0 & 0 & -2f(x_E(t),0) & 0
\end{array}
\right)
\left(
\begin{array}{c}
\xi_1 \\
\eta_1 \\
\xi_2 \\
\eta_2
\end{array}
\right).
\label{VE:y:py:zero:n:2}
\end{equation}
This system can be divided in two components: the \emph{tangential} one, with variables $(\xi_1,\eta_1)$; and the \emph{normal} one, normal variational equation, with variables $(\xi_2,\eta_2)$. Notice that these two systems of variational equations appear uncoupled, fact which simplifies its resolution.
Certainly, the tangential variational equation is
$\dot{\xi}_1=\eta_1, \ \dot{\eta}_1=0,$
whose solution is
\[
\xi_1(t)=\eta_1(0) t + \xi_1(0), \qquad \eta_1(t)=\eta_1(0).
\]
On the other hand, the normal variational equation always reduces to
\[
\dot{\xi}_2 =\eta_2, \qquad \dot{\eta}_2= -2f(x_E(t),0)\, \xi_2,
\]
or, equivalently,
\begin{equation}\label{sodv}
\ddot{\xi}_2 + 2f(x_E(t),0) \, \xi_2 =0, \qquad \eta_2=\dot{\xi}_2.
\end{equation}
If we denote by
\begin{equation}\label{fundamental}
\Phi(t)= \left(
\begin{array}{cc}
\phi_{11}(t) & \phi_{12}(t) \\
\phi_{21}(t) & \phi_{22}(t)
\end{array}
\right),
\end{equation}
the fundamental matrix solution of the normal variational equation satisfying that $\Phi(0)=\Identidad_2$, the identity matrix, then we have that
\[
\left(
\begin{array}{c}
\xi_2(t) \\ \eta_2(t)
\end{array} \right) =
\left(
\begin{array}{cc}
\phi_{11}(t) & \phi_{12}(t) \\
\phi_{21}(t) &\phi_{22}(t)
\end{array}
\right)
\left(
\begin{array}{c}
\xi_2(0) \\ \eta_2(0)
\end{array} \right).
\]
Therefore,
\[
\frac{\partial \xi_1(t)}{\partial \eta_1(0)}=t, \qquad \frac{\partial \xi_1(t)}{\partial \eta_2(0)}=0,
\qquad
\frac{\partial \xi_2(t)}{\partial \eta_1(0)} = 0, \qquad
\frac{\partial \xi_2(t)}{\partial \eta_2(0)} = \phi_{12}(t),
\]
and
\[
\xi_{2}(t)=c_1\xi_2^{(1)}(t)+c_2\xi_2^{(2)}(t),
\]
where $\{ \xi_{2}^{(1)},\xi_{2}^{(2)} \}$ is a basis of solutions of \eqref{sodv}.
So, from the formula~\eqref{def:J}, the Van Vleck-Morette matrix $J(t_1,t_0=0)$ reads
\begin{equation}%\label{J2}
J(t_1,0) =
\left(
\begin{array}{cc}
\frac{\partial \xi_1(t_1)}{\partial \eta_1(0)} & \frac{\partial \xi_1(t_1)}{\partial \eta_2(0)} \\[1.4ex]
\frac{\partial \xi_2(t_1)}{\partial \eta_1(0)} & \frac{\partial \xi_2(t_1)}{\partial \eta_2(0)}
\end{array}
\right)
=
\left(
\begin{array}{cc}
t_1 & 0 \\
0 & \phi_{12}(t_1)
\end{array}
\right),
\label{matrix:J}
\end{equation}
and has determinant $\det J(t_1, 0)=t_1\phi_{12}(t_1)$. From now on, we refer as variational equation the normal variational equation. Thus, the semiclassical approximation of the propagator along the solution $\gamma_E$, given by the Pauli-Morette formula \eqref{Pauli}, is
\begin{equation}
\label{paulin}
\WKB(x_1,t_1 \,|\, x_0,0) = \frac{1}{ 2\pi \rmi \hslash}\, \frac{1}{\sqrt{t_1\phi_{12}(t_1)}}\, {\rm e}^{\frac{\rmi}{2\hslash t_1} (x_1-x_0)^2 },
\end{equation}
As the path is on the line $y=0$ we will write, by abusing of notation,
$K_{\textsf{WKB}}(x_1,t_1 \,|\, x_0,0)$ to indicate
$K_{\textsf{WKB}}(\xbf_1,t_1 \,|\, \xbf_0,0)$.
\medskip

If $k>2$, straightforward computations lead to $\phi_{12}(t_1)=t_1$ and therefore $\det J(t_1, 0)=t_1^2.$ Hence,
\[
\WKB(x_1,t_1 \,|\, x,0) = \frac{1}{ 2\pi \rmi \hslash t_1}\,  {\rm e}^{\frac{\rmi}{2\hslash t_1} (x_1-x_0)^2 },
\]
which is the well known expression for the free particle propagator in two degrees of freedom.

\medskip

In the case $k=2$, the function $\phi_{12}(t_1)$ in~\eqref{matrix:J} can be determined from the values at $t_0=0$ and $t_1$ of a basis of solutions, $\{ \xi_2^{(1)}, \xi_2^{(2)} \}$. Precisely, in the general form
\[
\xi_2(t)=c_1 \xi_2^{(1)}(t) +  c_2 \xi_2^{(2)}(t), \qquad
\eta_2(t)=\dot{\xi}_2(t)=c_1 \dot{\xi}_2^{(1)}(t) +  c_2 \dot{\xi}_2^{(2)}(t),
\]
the values $c_1$ and $c_2$ are uniquely determined from the initial conditions at $t=0$ and so they are functions of $\xi_2(0)$ and $\eta_2(0)$ or, equivalently, of
$\xi_2^{(j)}(0)$ and $\dot{\xi}_2^{(j)}(0)$, for $j=1,2$.
So,
\begin{equation}
\phi_{12}(t)=\frac{\partial\xi_2(t)}{\partial \eta_2(0)} = \frac{\partial \xi_2(t)}{\partial c_1} \cdot
\frac{\partial c_1}{\partial \eta_2(0)} + \frac{\partial \xi_2(t)}{\partial c_2} \cdot
\frac{\partial c_2}{\partial \eta_2(0)} =
\xi_2^{(1)}(t) \cdot
\frac{\partial c_1}{\partial \eta_2(0)} + \xi_2^{(2)}(t) \cdot
\frac{\partial c_2}{\partial \eta_2(0)}.
\label{eq:dxi2deta20}
\end{equation}

To compute these two partial derivatives, we solve the linear system
\[
\xi_2(0)=c_1 \xi_2^{(1)}(0) +  c_2 \xi_2^{(2)}(0), \qquad
\eta_2(0)=c_1 \dot{\xi}_2^{(1)}(0) +  c_2 \dot{\xi}_2^{(2)}(0),
\]
by Cramer's rule (because it has a unique solution) and get
\[
c_1 = \frac{1}{D} \left( \xi_2(0) \dot{\xi}_2^{(2)}(0) - \eta_2(0) \xi_2^{(2)}(0) \right), \qquad \qquad
c_2 = \frac{1}{D} \left( \eta_2(0) \xi_2^{(1)}(0) - \xi_2(0) \dot{\xi}_2^{(1)}(0) \right),
\]
where
\[
D = \left|
\begin{array}{cc}
\xi_2^{(1)}(0) & \xi_2^{(2)}(0) \\[1.3ex]
\dot{\xi}_2^{(1)}(0) & \dot{\xi}_2^{(2)}(0)
\end{array}
\right|.
\]
Consequently,
\[
\frac{\partial c_1}{\partial \eta_2(0)} = - \frac {\xi_2^{(2)}(0)} D, \qquad \qquad
\frac{\partial c_2}{\partial \eta_2(0)} = \frac{\xi_2^{(1)}(0)}D,
\]
and substituting into expression~\eqref{eq:dxi2deta20} we obtain
\begin{equation*}
\phi_{12}(t_1) = \frac{\partial \xi_2(t_1)}{\partial \eta_2(0)} =  \frac 1 D
\left|
\begin{array}{cc}
\xi_2^{(1)}(0) & \xi_2^{(1)}(t_1) \\[1.3ex]
\xi_2^{(2)}(0) & \xi_2^{(2)}(t_1)
\end{array}
\right|.
%\label{eq:phi12:xi2}
\end{equation*}

Since the relevant variational equation in our study is the normal one (already denoted by, just, variational), from now we will remove from its basis of solutions the subscript $2$, that is, $\xi_2^{(j)}$
will be referred, simply, as $\xi^{(j)}$, $j=1,2$.
Thus, summarising, for $k=2$ our semiclassical approximate propagator $\WKB$ is given by the formula~\eqref{paulin}, where
\begin{equation}
\phi_{12}(t_1)=
\frac{1}{D}
\left|
\begin{array}{cc}
\xi^{(1)}(0) & \xi^{(1)}(t_1) \\[1.3ex]
\xi^{(2)}(0) & \xi^{(2)}(t_1)
\end{array}
\right|,
\qquad \qquad
D = \left|
\begin{array}{cc}
\xi^{(1)}(0) & \xi^{(2)}(0) \\[1.3ex]
\dot{\xi}^{(1)}(0) & \dot{\xi}^{(2)}(0)
\end{array}
\right|,
\label{eq:phi12:xi2}
\end{equation}
and $\{ \xi^{(1)}, \xi^{(2)} \}$ being a fundamental solution of the variational equation~\eqref{sodv}
\begin{equation*}
\ddot{\xi} + 2f(x_E(t),0) \, \xi =0,
\label{NVE:ode:2}
\end{equation*}
with $x_E(t)=\sqrt{2E} \, t + x_0$. Dot denotes derivative with respect to $t$.

This notation will be maintained henceforth in the paper.

%%%%%%%%%%%%%%%%%%%%%%%%%%%%%%%%%%%%%%%%%%%%%%%%%%%%%%%%%%%%%

\section{Applications}
\label{appl}
This section is devoted to the application of this result to some relevant families of equations. The type of the function $2f(x,0)$ appearing in the variational equation~\eqref{sodv} determines the family, according to the following table:
\begin{equation}\label{table}
\begin{array}{ccl}
2f(x,0)   &  \qquad &  \textrm{Family to which it reduces} \\
\hline \\[-0.5ex]
1-ax^2    &  \qquad & \quad \textrm{Hermite} \\[1.2ex]
b-\dfrac{a}{x^2}    & \qquad  & \quad \textrm{Bessel} \\[1.5ex]
-b + \dfrac{a}{\cosh^2 x}    & \qquad  & \quad \textrm{Legendre} \\[1.5ex]
-b - a \wp(x+\omega_3) & \qquad & \quad \textrm{Lam\'e}
\end{array}
\end{equation}
where in all four cases $a,b$ are real parameters. Here $\wp$ is the Weierstrass function  with real period $2\omega_1$
and imaginary period $2\omega_3$.

In the case $a=0$ all the potentials given in the table above have variational equations which reduce to a constant coefficients ode:
\[
\ddot{\xi} + \omega \xi =0,  \quad \omega \in \Rset.
\]
In particular the  $\WKB$ is integrable (in the sense of the Differential Galois theory) for any value of the energy $E$. For them, we have:
\arraycolsep=1.8pt
\def\arraystretch{2.5}
\[
\begin{array}{|c|c|c|c|c|} \hline
\quad \omega \quad & \{\xi^{(1)}, \xi^{(2)} \} & \phi_{12}(t_1) & \det J & \WKB(x_1,t_1|x_0,0) \\ \hline
0 & \{ 1,t\} & t_1 & t_1^2 & \dfrac{1}{2\pi\rmi \hslash} \, \dfrac{1}{|t_1|} e^ {\frac{\rmi}{2\hslash t_1} (x_1-x_0)^2 } \\[1.2ex] \hline
>0 & \{ \cos\sqrt{\omega}\, t, \sin\sqrt{\omega}\,t \} &
\dfrac{\sin\sqrt{\omega}\, t_1}{\sqrt{\omega}} & \dfrac{t_1\sin\sqrt{\omega}\, t_1}{\sqrt{\omega}} &
\dfrac{1}{2\pi\rmi \hslash} \, \dfrac{\omega^{1/4}}{\sqrt{t_1 \sin\sqrt{\omega}\, t_1}} e^{ \frac{\rmi}{2\hslash t_1} (x_1-x_0)^2} \\[1.2ex] \hline
<0 & \{ \cosh\sqrt{-\omega}\, t, \sinh\sqrt{-\omega}\,t \} &
\dfrac{\sinh\sqrt{-\omega}\, t_1}{\sqrt{-\omega}} & \dfrac{t_1\sinh\sqrt{-\omega}\, t_1}{\sqrt{-\omega}} &
\dfrac{1}{2\pi\rmi \hslash} \, \dfrac{(-\omega)^{1/4}}{\sqrt{t_1 \sinh\sqrt{-\omega}\, t_1}} e^{ \frac{\rmi}{2\hslash t_1} (x_1-x_0)^2} \\[1.2ex] \hline
\end{array}
\]

Considering $V(x,y)=\frac{\omega}{2}y^2$, the case $\omega=0$ is the free particle case and $\omega>0$ is the harmonic oscillator, which already is considered in the seminal works of Feynman, see for instance, \cite{fh2010,f2005}. Thus, the semiclassical approximation is the complete approximation: $\WKB=K$.

Henceforth in the paper, we will assume $a\neq 0$.

\subsection{Variational equations with Hermite equation. Verhulst potentials}

Potentials of the form $$V(x,y)=y^2 f(x,y), \quad f(x,y)=\frac{1}{2}(1-ax^2)+\hot (y),$$ where $\hot(y)$ means higher order terms in $y$. These systems are common in many physical systems. One example is the so-called
Verhulst's potentials, which take the form
\begin{equation}\label{verh}
V(x,y)=\frac{1}{2}\left(\omega_1^2x^2+\omega_2^2y^2\right)-\left(\frac{A_1}{3}x^3+A_2xy^2\right)-\left(\frac{B_1}{4}x^4+\frac{B_2}{2}x^2y^2+\frac{B_3}{4}y^4\right),
\end{equation}
with $\omega_1,\omega_2, A_1,A_2,B_1,B_2$, and $B_3$, real parameters. They were introduced by F. Verhulst~\cite{verhulst1979} to study systems of axi-symmetric galaxies. They exhibit a discrete-symmetric potential and can undergo resonances of type 1:2, 1:1, 2:1 and 1:3.
A suitable choice of the parameters and some trivial algebraic manipulations can lead them to fall into our class of potentials. For instance, swapping the variables $x$ and $y$,
\[
\tilde{V}(x,y)=V(y,x)=
\frac{1}{2}\left(\omega_1^2 y^2+\omega_2^2 x^2\right)-\left(\frac{A_1}{3}y^3+A_2yx^2\right)-\left(\frac{B_1}{4}y^4+\frac{B_2}{2}y^2x^2+\frac{B_3}{4}x^4\right),
\]
and taking $\omega_1=1$, $\omega_2=0$, $A_1=a_1$, $A_2=0$, $B_1=b_1$, $B_2=a$, and $B_3=0$, one gets
$\tilde{V}(x,y)=y^2 f(x,y)$, with
\begin{equation}
f(x,y)=\frac{1}{2}-\frac{a_1}{3}y-\frac{b_1}{4}y^2-\frac{a}{2}x^2 \qquad \textrm{and} \qquad f(x,0)=\frac{1}{2} - \frac{a}{2}x^2.
\label{verhulst:f}
\end{equation}
Notice that $\Gamma=\{ y=p_y=0 \}$ is an invariant plane of this system and hence the assumptions of Section~\ref{se:2dofhamiltonians} are satisfied.
The corresponding variational equation around $\Gamma$ is given by
\begin{equation}
\frac{d^2\xi}{dt^2} = - 2 f(x_E(t),0) \, \xi, \qquad \qquad \eta(t)=\frac{d\xi}{dt}(t),
\label{eq:nve:verhulst:1}
\end{equation}
with
\begin{eqnarray}
-2f(x_E(t),0) &=& -1 + a x_E^2(t) =
-1 +a \left( \sqrt{2E} \, t+x_0 \right)^2
\nonumber \\
&=&
2aE \, t^2 \pm 2\sqrt{2E} \, ax_0 \, t + \left( a x_0^2 -1 \right).
\label{VE:element:2}
\end{eqnarray}
The time transformation
\begin{equation}
s = \sqrt[4]{2Ea} \left( t \pm \frac{x_0}{\sqrt{2E}}  \right)
\label{prop:change:time}
\end{equation}
brings equation~\eqref{eq:nve:verhulst:1} with~\eqref{VE:element:2} into the harmonic oscillator ode
\begin{equation}
\frac{d^2 \xi}{ds^2}(s) = \left( s^2 - \lambda \right) \, \xi(s),
\qquad \qquad
\eta(s) = \sqrt[4]{2Ea} \, \frac{d\xi}{ds}(s),
\label{harmonic:oscil:ode}
\end{equation}
where $\lambda= \dfrac{1}{\sqrt{2Ea}}$. Observe that this transformation is equivalent to say
\[
s  = \pm \sqrt[4]{\frac{a}{2E}} \, x_E(t),
\]
and so $s$ is proportional to the spatial position of $x_E(t)$. This change simplifies the form of the variational equation. From the differential Galois Theory, it is know that ode~\eqref{harmonic:oscil:ode} admits Liouvillian solutions if and only if $\lambda=2m+1$, where $m\in \Nset \cup \{ 0 \}$. This implies that
the set of admissible ("Liouvillian", say) energies $E$ is discrete and it is given by
\begin{equation}
E_m = \frac{1}{2a} \left( \frac{1}{2m+1} \right)^2, \qquad  m\in \Nset \cup \left\{ 0 \right\}.
\label{admissible:E}
\end{equation}
Moreover, from relation~\eqref{energy:x0:x1} it follows that $(x_1-x_0)^2 = 2E_m t_1^2$ and so, fixed the initial position $x_0$, the unique admissible positions $x_1$ are those satisfying that
\begin{equation}
|x_1-x_0| = \frac{t_1}{\sqrt{a}} \frac{1}{2m+1},
\qquad m\in \Nset \cup \left\{ 0 \right\}.
\end{equation}
Regarding the solvability of~\eqref{harmonic:oscil:ode}, Galois Theory ensures their Liouvillian solutions to be of the form
\begin{equation}
\xi(s) = P_m(s) \, \rme^{-s^2/2},
\label{prop:xi2:Pm:m:pos}
\end{equation}
where $P_m(s)$ is a polynomial of degree $m$ (which we can assume, without loss of generality, to be monic). These polynomials $P_m$ satisfy the celebrated
Hermite differential equation
\begin{equation}
P_m''(s) - 2 s P_m'(s) + 2m P_m(s)=0.
\label{prop:ode:Pm}
\end{equation}
Its solutions $P_m$ are called
Hermite polynomials and, among other nice properties, they have the same parity as $m$. i.e.,
if $m$ is even then $P_m(s)$ is an even function and if $m$ is odd then $P_m(s)$ is an odd function (see~\cite{ismail2005}).

From D'Alembert formula, we know that
\[
\xi^{(2)}(s) = \xi^{(1)}(s) \, \int_0^{s} \frac{\rme^{z^2}}{P_m^2(z)}\, dz,
\]
is another solution of~\eqref{harmonic:oscil:ode} independent of $\xi^{(1)}$.
Together, $\{ \xi^{(1)}(s), \xi^{(2)}(s)\}$ form a fundamental solution.
In order to simplify the computations, we will assume henceforth in this section that $x_0=0$ and hence
\begin{equation}
s = \sqrt[4]{2E_m a} \, t.
\label{t:tau:hermite}
\end{equation}
The case $x_0\ne 0$ follows analogously.

In order to obtain an expression for $\phi_{12}(t_1)$ in $\det J(0,t_1)$ we apply formula~\eqref{eq:dxi2deta20}. In this context it reads as follows:
\begin{equation}
\phi_{12}(t_1) = \frac{1}{D}  \left|
\begin{array}{cc}
\xi^{(1)}(0) & \xi^{(1)}(s_1) \\[1.3ex]
\xi^{(2)}(0) & \xi^{(2)}(s_1)
\end{array}
\right|,
\label{phi12:hermite}
\end{equation}
where
\begin{equation}
D= \left|
\begin{array}{cc}
\xi^{(1)}(t=0) & \xi^{(2)}(t=0) \\[1.3ex]
\dot{\xi}^{(1)}(t=0) & \dot{\xi}^{(2)}(t=0)
\end{array}
\right| =
\frac{1}{\sqrt{2m+1}}
\left|
\begin{array}{cc}
\xi^{(1)}(s=0) & \xi^{(2)}(s=0) \\[1.3ex]
\dfrac{d\xi^{(1)}}{ds}(s=0) & \dfrac{d\xi^{(2)}}{ds}(s=0)
\end{array}
\right|,
\label{D:hermite}
\end{equation}
and we have taken into account
the relation between $s$ and $t$ given by~\eqref{t:tau:hermite}
and that
\begin{equation}
\sqrt[4]{2E_m a }=\frac{1}{\sqrt{2m+1}}, \qquad s= \frac{t}{\sqrt{2m+1}}.
\label{t:tau_x0:0}
\end{equation}
The parity of $P_m(s)$ determines (as it will be seen later), a separate study for $m$ even and $m$ odd.

\subsubsection{Case \boldmath{$m$} even}
Recall that in this case $P_m(s)$ is always an even function. Thus, on one side,
\[
\xi^{(1)}(s=0)=P_m(0) \ne 0, \qquad
\xi^{(2)}(s=0)=0,
\]
and on the other, since $P_m'(s)$ is odd,
\[
\frac{d\xi^{(1)}}{ds}(0) =0, \qquad  \frac{d\xi^{(2)}}{ds}(0) = \frac{1}{P_m(0)}.
\]
substituting in~\eqref{D:hermite} we get $D=\dfrac{1}{\sqrt{2m+1}}$. Moreover,
\[
\xi^{(1)}(s_1)= P_m(s_1) \, \rme^{-s_1^2/2}, \qquad
\xi^{(2)}(s_1) = P_m(s_1) \, \rme^{-s_1^2/2} \int_0^{s_1} \frac{\rme^{z^2}}{P_m^2(z)} \, dz,
\]
so in equation~\eqref{phi12:hermite} it gives rise to
\[
\phi_{12}(t_1) = \sqrt{2m+1}\, P_m(0) \, P_m(s_1) \, \rme^{-\frac{t_1^2}{2(2m+1)}}
\int_0^{\frac{t_1}{\sqrt{2m+1}}} \frac{\rme^{z^2}}{P_m^2(z)} \, dz.
\]
Consequently, the determinant of Van Vleck-Morette reads
\begin{equation}
\det J(0,t_1)= t_1 \phi_{12}(t_1) =
\sqrt{2m+1}\, P_m(0) \, P_m \left(\frac{t_1}{\sqrt{2m+1}} \right) \, \rme^{-\frac{t_1^2}{2(2m+1)}}
\int_0^{\frac{t_1}{\sqrt{2m+1}}} \frac{\rme^{z^2}}{P_m^2(z)} \, dz.
\label{hermite:detJ:m:even}
\end{equation}

\subsubsection{Case \boldmath{$m$} odd}

Same as in the previous case, we consider
${\ds \xi^{(1)}(s) = P_m(s) \rme^{-s^2/2}}$
as one solution of~\eqref{harmonic:oscil:ode}. The second one, provided by D'Alembert formula, is taken as
\begin{equation}
\xi^{(2)}(s) = \xi^{(1)}(s) \, \int_{s_0}^{s} \frac{\rme^{z^2}}{P_m^2(z)}\, dz,
\label{hermite:m:odd:xi2}
\end{equation}
with $s_0\ne 0$ to avoid the singularity inside the integral at $s=0$ (remind that since $P_m$ is odd then $P_m(0)=0$. Let us assume write
$P_m(s)=s^m + \cdots + a_3 s^3 + a_1 s$, where $a_1\ne 0$ for any odd $m$. Then, having in mind the Taylor expansion of the exponential function, we have that
\begin{equation}
\xi^{(1)}(s) =  a_1 s \left( 1 - \mathcal{O}(s^2) \right).
\label{prop:dem:m:odd:txi1}
\end{equation}
Moreover,
\[
\frac{\rme^{z^2}}{P_m^2(z)} = \frac{1}{a_1^2 z^2} \left( 1+ \left( 1-2\tilde{a}_3\right) z^2 + \mathcal{O}(z^4) \right),
\]
where $\tilde{a}_3 = \frac{a_3}{a_1}$,
and so
\[
\int_{s_0}^{s} \frac{\rme^{z^2}}{P_m^2(z)} \, dz =
\left( -\frac{1}{a_1^2 s} + \frac{1-2\tilde{a}_3}{a_1^2} s + \mathcal{O}(s^3) \right) - C_0,
\]
with
\[
C_0 = -\frac{1}{a_1^2 s_0} + \frac{1-2\tilde{a}_3}{a_1^2} s_0 +\cdots,
\]
is a constant obtained by evaluating any primitive function of
$\frac{\rme^{z^2}}{P_m^2(z)}$ at the point $s_0$.
Having in mind~\eqref{prop:dem:m:odd:txi1} and substituting the latter expression in~\eqref{hermite:m:odd:xi2} it follows that
\begin{eqnarray*}
\xi^{(2)}(s) &=& \left( a_1 s + a_3 s^3 + \mathcal{O}(s^5) \right) \cdot \left( 1 - \frac{s^2}{2} +
\mathcal{O}(s^4) \right) \cdot \left(
\left( -\frac{1}{a_1^2 s} + \frac{1-2\tilde{a}_3}{a_1^2} s + \mathcal{O}(s^3) \right) - C_0 \right)
\\
&=& -\frac{1}{a_1} - a_1 C_0 s + \left( -\frac{a_3}{a_1^2} + \frac{3}{2a_1}
- \frac{2 \tilde{a}_3}{a_1} \right) s^2 + \mathcal{O}(s^3).
\end{eqnarray*}
Since $a_1 \ne 0$, the limit at $s=0$ is well defined. Indeed,
$\lim_{s \rightarrow 0} \xi^{(2)}(s) = -\frac{1}{a_1}$.
This implies that $\xi^{(2)}(s)$ admits an analytic extension in $\Cset$. Abusing of notation, we denote this extension with the same name and so we write
\[
\xi^{(2)}(0) =-\frac{1}{a_1}.
\]
We have now the ingredients to compute $\phi_{12}(t_1)$ through formulas~\eqref{phi12:hermite}-\eqref{D:hermite}. Indeed, evaluating at $s=0$,
\[
\xi^{(1)}(0)=P_m(0)=0, \qquad \xi^{(2)}(0)=-\frac{1}{a_1}, \qquad \frac{d\xi^{(1)}}{ds}(0)=P_m'(0)=a_1, \qquad \frac{d\xi^{(2)}}{ds}(0) = -a_1 C_0
\]
it follows that ${\ds D= \frac{1}{\sqrt{2m+1}}. }$ Therefore,
\[
\phi_{12}(t_1) = \sqrt{2m+1} \left|
\begin{array}{cc}
\phantom{+}0  & \xi^{(1)}(s_1) \\[1.2ex]
-\dfrac{1}{a_1} & \xi^{(2)}(s_1)
\end{array}
\right| = \frac{\sqrt{2m+1}}{P_m'(0)} \, P_m \left( \frac{t_1}{\sqrt{2m+1}} \right)  \, \rme^{- \frac{t_1^2}{2(2m+1)}},
\]
and so
\begin{equation}
\det J(0,t_1) =  \frac{\sqrt{2m+1}}{P_m'(0)} \, P_m \left( \frac{t_1}{\sqrt{2m+1}} \right)  \, t_1 \, \rme^{- \frac{t_1^2}{2(2m+1)}}.
\label{hermite:detJ:m:odd}
\end{equation}

%\subsubsection{An example}
%\label{verhulst:example}
To illustrate this case we propose the following simple example. \medskip

\noindent{\textbf{Example.}}
Let us assume our potential to be $V(x,y)=y^2 f(x,y)$ with
\[
f(x,y)=\frac{1}{2}-\frac{b_1}{4}y^2-\frac{a}{2}x^2.
\]
The Hamiltonian function is
\[
H(x,p_x,y,p_y) = \frac{p_x^2+p_y^2}{2} + y^2 f(x,y),
\]
and its associated (Verhulst) Hamiltonian system
\[
\begin{array}{lclclcl}
\dot{x}   &=& p_x   & \qquad & \dot{y} &=& p_y  \\
\dot{p}_x &=& xy^2  & \qquad & \dot{p}_y &=& - y \left( 1 - b_1 y^2 - a x^2 \right)
\end{array}
\]
Let us fix an energy level $H=E>0$ and consider a free particle moving on the invariant plane $\Gamma=\{ y=p_y=0\}$ between two different points $(x_0,p_{x_0},0,0)$ (at time $t_0=0$) and $(x_1,p_{x_1},0,0)$ (at time $t=t_1$). Its trajectory is given by $(x_E(t),p_E(t),0,0)$ where $x_E(t)=\sqrt{2E} \,t$ and $p_E(t)=\dot{x}_E(t)=\sqrt{2E}$. Therefore, there exists a discrete set of values of the energy
\[
\left\{ E=E_m= \frac{1}{2} \left( \frac{1}{2m+1} \right)^2 \ \bigg| \ m \in \Nset \cup \{0\} \right\},
\]
for which the determinant of Van Vleck-Morette admits a closed expression in terms of Liouvillian functions. This determinant is given by~\eqref{hermite:detJ:m:even} if $m$ even and by~\eqref{hermite:detJ:m:odd} if $m$ odd. In the particular case that $m=0$, that is $E_0=1/2$, after simplification due to $t_1=|x_1-x_0|$, we obtain
\[
\WKB(x_1,t_1 \,|\, x_0,0) = \frac{1}{ 2\pi^{5/4} \rmi \hslash}
\, \frac{\rme^{t_1^2/4}}{\sqrt{t_1 \erfi(t_1)}}
\, \rme^{\frac{\rmi}{2\hslash} \, |x_1-x_0| },
\]
where
\[
\erfi(t)=-\rmi \, \erf(\rmi\, t_1) = \frac{2}{\sqrt{\pi}} \int_0^t \rme^{s^2} \, ds,
\]
is the imaginary error function.

%%%%%%%%%%%%%%%%%%%%%%%%%%%%%%%%%%%%%%%%%%%%%%%%%%%%%%%%%%%%%

\subsection{Variational equations with Bessel functions}

Rational potentials have been widely studied by many authors (see, for example,~\cite{ab2008} and references therein). They can fall in the hypothesis of this work by considering, for instance, potentials $V$ of the form
\begin{equation}\label{rmbesve}
V(x,y)=y^2f(x,y),\quad f(x,y)=\frac{b}{2}-\frac{a}{2x^2}+ \hot (y),
\end{equation}
$a,b$ being real parameters non-simultaneously vanishing. This kind of potentials have singularities, which in our case, correspond to $x=0$. To avoid it, in our study we will assume that $0 \notin [x_0,x_1]$, the $x$-interval where our free particle moves.

By the using of the change of variable $\tau=x_E(t)=\sqrt{2E}t+x_0$ in the variational equation \eqref{sodv}  we obtain
\begin{equation}\label{sodvalg}
\frac{d^2\xi_2}{d\tau^2}=\frac{1}{2E}r(\tau)\xi_2, \quad r(\tau)=-2f(\tau,0).
\end{equation}

Following the procedure given in Section~\ref{se:2dofhamiltonians}, for a given value $E>0$ of the energy and considering the free particle motion given by~\eqref{freeparticle:xH},
we reduce the study to its variational equation~\eqref{sodvalg}, which now reads
\begin{equation}\label{nve:bessel0}
\frac{\partial^2\xi}{\partial \tau^2}=\frac{1}{2E}\left(\frac{a}{\tau^2}-b\right)\xi.
\end{equation}
Equation \eqref{nve:bessel0}, for $b\neq 0$, is one of the equivalent expressions of the well-known Bessel differential equation. For further details about Bessel odes see, for instance, Appendix~\ref{apbes} and the classical reference~\cite{watson1922}.

Comparing~\eqref{nve:bessel0} with the standard Bessel equation in normal form
\begin{equation}\label{nve:stanbessel}
\frac{\partial^2\xi}{\partial \tau^2}=\left( \frac{\nu(\nu+1)}{\tau^2}- \mu^2 \right) \xi,\,\qquad  \nu \in \Cset, \qquad \mu \neq 0,
\end{equation}we get
\[
a=2\nu(\nu+1)E, \qquad \qquad b=2E \mu^2.
\]
It is well known that Bessel equation \eqref{nve:stanbessel} is integrable if and only if $\nu:=n$ is an integer number (see Appendix A). In this case, the set of values of the energy for which we have Liouvillian solutions is given by
\[
\left\{ E_n = \frac{|a|}{2n(n+1)} \ \bigg| \ n \in \Zset\setminus \{ -1,0\} \right\}.
\]
Thus, for $\nu=n$, a basis of solutions of~\eqref{nve:stanbessel} is given by $\xi^{(1)}(t_1)=\sqrt{x_1}J_{n+\frac{1}{2}}(x_1)$ and $\xi^{(2)}(t_1)=\sqrt{x_1}Y_{n+\frac{1}{2}}(x_1)$. Then by equation \eqref{eq:phi12:xi2}  we have $\phi_{12}(t_1)$ as follows:
$$\phi_{12}(t_1)=\frac{t_1}{x_1-x_0}\sqrt{\frac{x_1}{x_0}}\left(\frac{J_{n+\frac{1}{2}}(x_0)Y_{n+\frac{1}{2}}(x_1)-J_{n+\frac{1}{2}}(x_1)Y_{n+\frac{1}{2}}(x_0)}{J_{n+\frac{1}{2}}(x_0)Y_{n-\frac{1}{2}}(x_0)-J_{n-\frac{1}{2}}(x_0)Y_{n+\frac{1}{2}}(x_0)}\right).$$

Thus, the  Van Vleck-Morette determinant reads $$\det J(t_1,0)=t_1\phi_{12}(t_1)=\frac{t_1^2}{x_1-x_0}\sqrt{\frac{x_1}{x_0}}\left(\frac{J_{n+\frac{1}{2}}(x_0)Y_{n+\frac{1}{2}}(x_1)-J_{n+\frac{1}{2}}(x_1)Y_{n+\frac{1}{2}}(x_0)}{J_{n+\frac{1}{2}}(x_0)Y_{n-\frac{1}{2}}(x_0)-J_{n-\frac{1}{2}}(x_0)Y_{n+\frac{1}{2}}(x_0)}\right).$$

Finally, the semiclassical approximation for any value of $x_1,x_0,t_1$ of the propagator is given by Pauli-Morette formula:
\[
K_{\textsf{WKB}}({\bf x}_1,t_1 \,|\, {\bf x}_0,0)=\frac{1}{ 2\pi \rmi \hslash}\, \frac{\sqrt{x_1-x_0}}{t_1}\sqrt[4]{\frac{x_0}{x_1}}\left(\frac{J_{n+\frac{1}{2}}(x_0)Y_{n-\frac{1}{2}}(x_0)-J_{n-\frac{1}{2}}(x_0)Y_{n+\frac{1}{2}}(x_0)}{J_{n+\frac{1}{2}}(x_0)Y_{n+\frac{1}{2}}(x_1)-J_{n+\frac{1}{2}}(x_1)Y_{n+\frac{1}{2}}(x_0)}\right) \rme^{\frac{\rmi}{\hslash} \frac{(x_1-x_0)^2}{2t_1}}.
\]

%%%%%%%%%%%%%%%%%%%%%
\medskip

\textbf{Example.} To illustrate the previous case, we consider the Hamiltonian with potential \eqref{rmbesve} being $\nu=2$. The expressions for $J_{\frac{5}{2}}$, $Y_{\frac{5}{2}}$, (and so the ones for $J_{\frac{3}{2}}$ and $Y_{\frac{3}{2}}$) read as follows:
%\begin{center}
$$
\begin{array}{rclcrcl}
J_{\frac{5}{2}}(\tau) &=&{\ds -\sqrt{\frac{2}{\pi \tau}}\left(\sin(\tau)+\frac{3\cos(\tau)}{\tau} -\frac{3\sin(\tau)}{\tau^2}\right)}, &\qquad&
J_{\frac{3}{2}}(\tau) &=&{\ds \sqrt{\frac{2}{\pi \tau}}\left(\frac{\sin(\tau)}{\tau}-\cos(\tau)\right)}, \\[1.5ex]
Y_{\frac{5}{2}}(\tau) &=& {\ds \sqrt{\frac{2}{\pi \tau}}\left(\cos(\tau)-\frac{3\sin(\tau)}{\tau} -\frac{3\cos(\tau)}{\tau^2}\right)}, &\qquad& Y_{\frac{3}{2}}(\tau)&=&{\ds -\sqrt{\frac{2}{\pi \tau}}\left(\frac{\cos(\tau)}{\tau}+\sin(\tau)\right)}.
\end{array}
$$
%\end{center}
Thus, the semiclassical approximation of the propagator is given by Pauli-Morette- formula \eqref{Paulia}
$$
K_{\textsf{WKB}}( x_1,t_1 \,|\, x_0,0)=\dfrac{1}{ 2\pi \rmi \hslash}\, \frac{\sqrt{x_1-x_0}}{t_1}\sqrt[4]{\frac{x_0}{x_1}}\left(\frac{J_{\frac{5}{2}}(x_0)Y_{\frac{3}{2}}(x_0)-J_{\frac{3}{2}}(x_0)Y_{\frac{5}{2}}(x_0)}{J_{\frac{5}{2}}(x_0)Y_{\frac{5}{2}}(x_1)-J_{\frac{5}{2}}(x_1)Y_{\frac{5}{2}}(x_0)}\right) \rme^{\frac{\rmi}{\hslash} \frac{(x_1-x_0)^2}{2t_1}}.
$$
Having in mind that $$J_{\frac{5}{2}}(x_0)Y_{\frac{3}{2}}(x_0)-J_{\frac{3}{2}}(x_0)Y_{\frac{5}{2}}(x_0)=\frac{2}{\pi x_0},$$ and setting $$G:=J_{\frac{5}{2}}(x_0)Y_{\frac{5}{2}}(x_1)-J_{\frac{5}{2}}(x_1)Y_{\frac{5}{2}}(x_0),$$
for $0<x_0<x_1$, we have
$$G=\frac{(2x_1^2 (x_0^2 - 3) + 18 x_1 x_0 - 6 x_0^2 + 18) \sin(x_1 - x_0) - 6 (x_1 - x_0) (x_1 x_0 + 3) \cos(x_1 - x_0)}{\pi (x_1x_0)^\frac{5}{2}}.
$$
Thus,
$$
K_{\textsf{WKB}}(x_1,t_1 \,|\, x_0,0)=\dfrac{1}{ \pi^2 \rmi \hslash}\, \frac{\sqrt{x_1-x_0}}{x_0t_1}\sqrt[4]{\frac{x_0}{x_1}}\frac{1}{G} \rme^{\frac{\rmi}{\hslash} \frac{(x_1-x_0)^2}{2t_1}}.
$$

%%%%%%%%%%%%%%%%%%%%%%%%%%%%%
\bigskip
%\vskip1cm
On the other side, if we set $b=0$ and $a=2\nu(\nu+1)E$ in
equation~\eqref{nve:bessel0}, with $\nu\in \Rset \setminus \{ 0,-1 \}$, we obtain an Euler-Cauchy differential equation
\begin{equation}\label{nve:cauchy}
\frac{\partial^2\xi}{\partial \tau^2}=\frac{\nu(\nu+1)}{\tau^2}\, \xi.
\end{equation}
Its general solution, for $\nu\neq -1/2$, is
$$\xi(\tau)= c_1\xi^{(1)}(\tau)+ c_2\xi^{(2)}(\tau)=c_1\tau^{\nu+1}+c_2\tau^{-\nu},\quad \mbox{with} \quad \xi^{(1)}(\tau)=\tau^{\nu+1},\quad  \xi^{(2)}(\tau)=\tau^{-\nu}.$$
As in the previous case, the expression for $\phi_{12}(t_1)$ is derived from equation~\eqref{eq:phi12:xi2}:
\[
\phi_{12}(t_1)=
\frac{x_1^{\nu+1}x_0^{-\nu}-x_0^{\nu+1}x_1^{-\nu}}{(2\nu+1)\sqrt{2E}}=\frac{x_0t_1}{(2\nu+1)(x_1-x_0)}\left(\left(\frac{x_1}{x_0}\right)^{\nu+1}-\left(\frac{x_1}{x_0}\right)^{-\nu}\right).
\]
The expression for the propagator $\WKB$ follows analogously:
$$
K_{\textsf{WKB}}( x_1,t_1 \,|\, x_0,t_0) = \frac{\sqrt{2\nu+1}}{ 2\pi \rmi \hslash}\, \frac{\sqrt{\frac{x_1}{x_0}-1}}{t_1\sqrt{\left(\frac{x_1}{x_0}\right)^{\nu+1}-\left(\frac{x_1}{x_0}\right)^{-\nu}}}\, \rme^{\frac{\rmi}{\hslash} \frac{(x_1-x_0)^2}{2t_1}}.
$$

\medskip

In the case $\nu=-1/2$ one has
\begin{comment}
$$
\xi_{2}^{(1)}(t_1)=\sqrt{x_1},\quad\xi_{2}^{(2)}(t_1)=\sqrt{x_1}\log\left(x_1\right), \quad \eta_{2}^{(1)}(t_1)=\frac{1}{2}\sqrt{\frac{2E}{x_1}}, \quad \eta_{2}^{(2)}(t_1)=\frac{1}{2}\sqrt{\frac{2E}{x_1}}\left(2+\log(x_1)\right).$$
Then by equation \eqref{eq:phi12:xi2}  we have:
\end{comment}
$$
\phi_{12}(t_1)=
\frac{\sqrt{x_0x_1}}{\sqrt{2E}}\left(\log(x_1)-\log(x_0)\right)=\frac{\sqrt{x_0x_1}t_1}{x_1-x_0}\log\left(\frac{x_1}{x_0}\right).$$
and the semiclassical approximation of the propagator reads
$$
K_{\textsf{WKB}}({ x}_1,t_1 \,|\, { x}_0,t_0) = \frac{1}{ 2\pi \rmi \hslash}\, \frac{\sqrt{
x_1-x_0}}{t_1\sqrt[4]{x_0x_1}{\sqrt{\log(x_1)-\log(x_0)}}}\, \rme^{\frac{\rmi}{\hslash} \frac{(x_1-x_0)^2}{2t_1}}.
$$
%%%%%%%%%%%%%%%%%%%%%%%%%%%%%%%%%%%

%%%%%%%%%%%%%%%%%%%%%%%%%%%%%%%%%%%%%%%%%%%%%%%%%%%%%%%%%%%%%

\subsection{Variational equations with Legendre functions}
The potentials considered in this section are those of the form \begin{equation}\label{rmpotve0}
V(x,y)=y^2f(x,y)  \qquad f(x,y)=-\frac{b}{2}+\frac{a}{2\cosh^2(x)}+\hot(y).
\end{equation}

Through the change of variable $\tau=x_E(t)=\sqrt{2E}t+x_0$ in~\eqref{sodv}, such as in the case of Bessel family,  we obtain the variational equation in the form~\eqref{sodvalg}.

For a given value of the energy $E>0$ and consider the free particle motion given by~\eqref{freeparticle:xH}.
This reduces the study to the one of its variational equation~\eqref{sodvalg}, which now reads
\begin{equation}\label{ex3}
\frac{\partial^2{\xi}}{\partial \tau^2}=\frac{1}{2E}\left(b-\frac{a}{\cosh^2(\tau)}\right)\xi.
\end{equation}
Comparing~\eqref{ex3} with the standard differential equation involving the so-called Rosen--Morse potential
\begin{equation}\label{ex3b}
\frac{\partial^2{\xi}}{\partial \tau^2}=\left(\mu^2-\frac{\nu(\nu+1)}{\cosh^2(\tau)}\right)\xi,
\end{equation}
we get the conditions
\[
\nu(\nu+1)=\frac{a}{2E},\qquad \mu^2=\frac{b}{2E},
\]
which determine the possible values of the energy.
Recall that the Rosen--Morse potentials were studied in~\cite{rosen1932} (see also~\cite{morales2015,aas2018}).
To apply differential Galois techniques we introduce the change of variable $z=\tanh(\tau)$ (see also~\cite{aas2018}) and obtain
\begin{equation}\label{LE3}
(1-z^2)\frac{\partial^2 \xi}{\partial z^2}-2z\frac{\partial \xi}{\partial z}+\left(\nu(\nu+1)-\frac{\mu^2}{1-z^2}\right)\xi=0.
\end{equation}
Equation~\eqref{LE3} is the well known Legendre equation.
Its integrability, in terms of their parameters, is given in Proposition~\ref{propLegendre} (see the Appendix~\ref{A2}) and it is only achieved for a discrete set of values of the parameters.

To illustrate the computation of the Feynman propagator for the Legendre family we restrict ourselves to the case $\mu=0$ (i.e. $b=0$).
The case $b\ne 0$ would follow a similar procedure, using also Proposition~\ref{propLegendre}.

Having in mind Proposition~\ref{propLegendre}(a), the case $\mu=0$  corresponds to consider $\nu=n\in\Zset$. Moreover, without loss of generality, we can assume $n\in \Nset$. Then, for
\[
E_n=\frac{a}{2n(n+1)}, \qquad \mu_n^2=\frac{b}{2E_n} \qquad n \in \Nset,
\]
the general solution of Legendre equation~\eqref{LE3} is given by
$$
\xi(z)=c_1\, P_n(z)+c_2\, Q_n(z),
$$
where $P_n(z)$ denotes the Legendre polynomial of degree $n$ and $Q_n(z)$ an independent solution determined by D'Alembert formula. Therefore, its general solution becomes
$$
\xi(\tau)=c_1\, P_n(\tanh(\tau))+c_2 \,Q_n(\tanh(\tau)), \qquad \tau=\frac{x_1 - x_0}{t_1}\, t+x_0.
$$
That is, $\{ \xi^{(1)}(\tau), \xi^{(2)}(\tau) \}= \{ P_n(\tanh(\tau)), \ Q_n(\tanh(\tau)) \}$ is a basis of solutions. Thus,
$$
\xi^{(1)}(0)=P_n(\tanh(x_0)),\quad \xi^{(2)}(0)=Q_n(\tanh(x_0)),\quad \xi^{(1)}(t_1)=P_n(\tanh(x_1)),\quad \xi^{(2)}(t_1)=Q_n(\tanh(x_1))
$$
and
$$
\dot\xi^{(1)}(t=0)=\frac{x_1-x_0}{t_1\cosh^2(x_0)}P_n'(\tanh(x_0)),\qquad \dot\xi^{(2)}(t=0)=\frac{x_1-x_0}{t_1\cosh^2(x_0)}Q_n'(\tanh(x_0)),
$$
where $'$ denotes, in this case, derivative with respect to $\tau$.
Applying formula~\eqref{eq:phi12:xi2} we obtain
$$
\phi_{12}(t_1)=\frac{t_1\cosh^2(x_0)}{x_1-x_0}\cdot\frac{P_n(\tanh(x_0))\,Q_n(\tanh(x_1))-P_n(\tanh(x_1))\,Q_n(\tanh(x_0))}{P_n(\tanh(x_0))\,Q'_n(\tanh(x_0))-P'_n(\tanh(x_0))\,Q_n(\tanh(x_0))},
$$
and the corresponding semiclassical approximation of the propagator is
$$
\WKB(x_1,t_1 \,|\, x_0,0) = \frac{1}{ 2\pi \rmi \hslash}\, \frac{1}{\sqrt{t_1\phi_{12}(t_1)}}\, {\rm e}^{\frac{\rmi}{2\hslash t_1} (x_1-x_0)^2 }.
$$

\subsection{Variational equations with Lam\'e functions}
This family corresponds to potentials of type
\begin{equation}\label{potLam}
V(x,y)=y^2f(x,y) \qquad
f(x,y)=\left(-\frac{b}{2}-\frac{a}{2}\wp(x+\omega_3)\right)+ \hot(y),
\end{equation}
where, on the free particle solution, it simply reads
\[
f(x_E(t),0)=-\frac{b}{2}-\frac{a}{2} \wp(x_E(t)+\omega_3).
\]
with $x_E(t)=\sqrt{2E}\,t+x_0=\frac{x_1-x_0}{t_1}t+x_0$.
Notice that the potential~\eqref{potLam} depends on four parameters: $a$, $b$, $g_2$ and $g_3$; the parameter $\omega_3$ depending, on its turn, on $g_2$ and $g_3$,
see Appendix \ref{apLame} and for more details see \cite{ww2020}.

By the change $\tau=x_E(t)+\omega_3$, the  variational equation~\eqref{sodv} becomes the celebrated Lam\'e differential equation
\begin{equation}\label{ex4}
\frac{d^2\xi}{d\tau^2} =\left(n(n+1)\wp(\tau)+B)\right)\xi,
\end{equation}
with
\[
n(n+1)=\frac a{2E}, \qquad B=\frac b{2E}.
\]

%%%%%%%%%%%%%%%
%%%%%%%%%%%%%%%%%
According to the Appendix\ref{apLame},~  the three cases of differential Galois integrability of Lam\'e equation lead only to a discrete set of values of parameters. Indeed,
\[
E=E_n=\frac{a}{2n(n+1)}, \qquad B=B_n=n(n+1)\frac{b}{a}.
\]
As we did for the Legendre family of potentials, we illustrate the theory by considering a concrete example. Namely, let us take $n=1$ and $a=2$. Hence, we our potential becomes
\begin{equation}\label{cpotLam}
V(x,y)=y^2f(x,y) \qquad
f(x,y)=\left(-\frac{b}{2}-\wp(x+\omega_3)\right)+ \hot(y).
\end{equation}
We recall that $\wp$ is the Weierstrass function, which is a solution of  the differential equation $(dw/dz)^2=h(w)$, with $h(w)=4w^3-g_2w-g_3$ and discriminant $\Delta=g_2^3-27g_3^2\neq 0$ in order that the polynomial $h(w)$ has simple roots (otherwise it can be transformed into simpler forms).

Here, we consider $h(B)\neq 0$, which falls into the Hermite-Halphen family, an integrable case of the Lam\'e equation (see Appendix~\ref{apLame} and also~\cite{ww2020,ms1996}).

%%%%%%%%%%%%%%%%%%%%%%%%
{\bf Example.} Fix $b=1$ in the potential~\eqref{cpotLam} and $h(w)=4w^3 - 28w + 24$, with roots $e_1=2$, $e_2=1$, and $e_3=-3$. The parameters are $g_2=28$, $ g_3=-24$ and (according to the notation of Appendix~\ref{apLame}) the discriminant is $\Delta=6400$. To ease the computations we fix $x_0=0$.
We have $E_1=1/2$ and define $\tau=t+\omega_3$. Moreover, $h(B)\neq 0$, i.e., $B\neq 1,2,-3$.
%Assume also $n=1$ and
Then equation~\eqref{ex4} becomes
\begin{equation}\label{ex6}
\ddot{\xi}=\left(B+2\wp(t)\right)\xi,
\end{equation}
with a basis of solutions $\{ \xi^{(1)}(t),\xi^{(2)}(t) \}$, where
$$
\xi^{(1)}(t)=\sqrt{B-\wp(t+\omega_3)}\,e^{\frac{1}{2}\sqrt{h(B)}\int_{0}^{t}\frac{ds}{B-\wp(s+\omega_3)}},\quad\xi^{(2)}(t)=\sqrt{B-\wp(t+\omega_3)}\,e^{-\frac{1}{2}\sqrt{h(B)}\int_{0}^{t}\frac{ds}{B-\wp(s+\omega_3)}}.
$$
For $t=0$, we have $\tau=\omega_3$ and therefore
$\xi^{(1)}(0)=\sqrt{B+3}$, $\xi^{(2)}(0)=\sqrt{B+3}.$
Thus, the value of $D$ in formula~\eqref{eq:phi12:xi2} is $D=2\sqrt{h(B)}$ and  $$\phi_{12}(t_1)=\sqrt{\frac{(B+3)\left(B-\wp(t_1+\omega_3)\right)}{h(B)}}\,\sinh\left(\frac{1}{2}\sqrt{h(B)}\int_0^{t_1}\frac{ds}{B-\wp(s+\omega_3)}\right).$$
Finally,
$$
\WKB(x_1,t_1 \,|\, x_0,0) = \frac{1}{ 2\pi \rmi \hslash}\, \frac{1}{\sqrt{t_1\phi_{12}(t_1)}}\, {\rm e}^{\frac{\rmi}{2\hslash t_1} (x_1-x_0)^2 }.
$$

%%%%%%%%%%%%%%%%%%%%%

\section{Final remarks: non-integrability and future work}\label{finrem}
%%%%%%%%%%%%%%%%%%%%%%%%%%%%%%%%%%%%%%%%%%%%%%%%%%%%%%%%%%%%%

In this section we precise the non-integrability statements concerning the families of potentials and the corresponding propagators considered in this paper and we present some questions for future work.

In fact, the families considering  here  are generically non-integrable.

\begin{proposition}\label{prop:noint}
Consider the families of potentials defined by equation~\eqref{def:potential} with $k=2$ and given in Table~\eqref{table}.

Then, the corresponding Hamiltonian systems whose variational equations reduce to Bessel and Legendre families with $b\neq 0$ and those which reduce to the
Hermite and Lam\'e families are not integrable in the Liouville - Arnold sense.
\end{proposition}

\pproof{The main common argument is the following: according to Section \ref{appl}, the variational equations of the Hamiltonian systems are differential Galois integrable only for discrete sets of the energy. Consequently, the Hamiltonian systems cannot be integrable in the Liouville -Arnold sense.

First, consider the Hermite and Bessel families. As there are values of the energy $E$ for which the corresponding variational equations are not integrable in the sense of differential Galois theory, then the Hermite and Bessel families of Hamiltonian systems are not integrable in the Liouville - Arnold sense (see~\cite{moralesramis2001, morales2013}).

Now we consider Legendre family. From Proposition~\ref{propLegendre} (see Appendix~\ref{A2}) we notice that only discrete values of the energy are compatible with integrability. Therefore the main argument above applies.

And last, but not least, we consider the Lam\'e family. For the case of Lam\'e and the Hermite-Halphen solutions, the Brioschi-Halphen-Crawford solutions as well as the Baldassarri solutions, they are only compatible with integrability for discrete values of the energy (for more details see~\cite[\S 3]{ms1996} and also~\cite{morales2013}). Analogously, the claim follows.
}

We stress that, in Proposition~\ref{prop:noint}, since for the Hermite and Bessel families the variational equations have irregular singular points at infinity, then the obstruction is to the existence of an additional rational first integral. For the other two families, the obstruction is to the existence of an additional meromorphic first integral (see~\cite{moralesramis2001,morales2013}).
In the philosophy of the papers \cite{m2020,as2013}, we introduce the following definition:
\begin{definition}\label{proint}
The semiclassical approximation of the Feynman propagator is integrable if for any fixed classical path $\gamma$ the expression $\WKB=\WKB(t_1)$ is a Liouvillian function over the coefficient field of the associated variational equation.
\end{definition}

In the above definition, the points $x_0$ and $x_1$ are considered fixed and so, the semiclassical aproximation of the Feynman propagator \eqref{Paulia} becomes a function depending only on time $t_1$.

If there exists a classical path whose $\WKB=\WKB(t)$ is not a Liouvillian function, then the semiclassical approximation of the Feynman propagator is not integrable.

Recall that a Liouvillian function over a differential field of functions $\mathcal{K}$ is a function obtained by a combination of algebraic functions, integrals and exponential of integrals of functions in $\mathcal{K}$. Then, the solutions of a linear differential equation over $\mathcal{K}$ is given by Liouvillian functions if, and only if,  the linear differential equation is integrable (for a more formal statement  see the appendix in~\cite{m2020}). The differential coefficient field $\mathcal{K}$ along this paper is generated by the function $f(x_E(t),0)$: see formula~\eqref{VE:y:py:zero:n:2}.  It is clear from
formula~\eqref{paulin} that $\WKB(t_1)$ is a Liouvillian function if, and only if, $\phi_{12}(t_1)$ is, on its turn, a Liouvillian function. But, by~\eqref{eq:phi12:xi2}, $\phi_{12}(t_1)$ is given as a linear combination of the base of solutions of the normal variational equation. Consequently, all the closed form formulas  of  $\phi_{12}(t_1)$ obtained in this paper are given by Liouvillian functions.

However, from the proof of the above proposition, it becomes clear that it will not be possible to obtain for some of the paths a base of solutions of the normal variational equation as Liouvillian functions. So, the following corollary can be stated:

\begin{corollary}\label{GOAL}
Consider the families of potentials defined by equation~\eqref{def:potential} with $k=2$ given in Table~\eqref{table}. Then the semiclassical approximations of the Feynman propagator of  Bessel and Legendre families with $b\neq 0$, as well as Hermite and Lam\'e families, are not integrable.
\end{corollary}
%%%%%%%%%%%%%%%%%%%%%%%%%%%%%%%%%%%%%%%%%%%%%%%%%%%%%%%%%%%%%

In our opinion, although this paper could be considered as an academic exercise, we believe that it can play a useful r\^ole as a first illustration of the applicability of the Differential Galois methodology in Quantum Physics. In this work we have used it to compute the semiclassical quantification around a simple type of classical paths.

%avoiding many technical complications of other more interesting physical problems.
%(no acabo de entender a qué nos referimos en esta frase).T

It seems to us that there are other, less academic, more challenging applications to physical problems. For instance, the ones related to the quantification of periodic orbits or those regarding \emph{tunneling}. More precisely, three of the problems we would like to address in future works are:

\begin{enumerate}
\item The quantification around families of periodic orbits parametrized by the energy, following essentially the work of Gutzwiller~\cite{Gutzwiller1971}.

\item The macroscopic tunneling problems in magnetism, as proposed in the book of Chudnovsky and Tejada~\cite{Chudnovsky-Tejada}. However most of the systems appearing therein are $1$-degree of freedom, the use of a non-standard Hamiltonian function (that is, not of the form kinetic + potential) can give rise to interesting dynamical features.

\item Tunneling problems in higher degrees of freedom systems, but with enough number of symmetries.

\end{enumerate}

\bigskip

\subsection*{Acknowledgements.}
PA-H is partially supported by the Ministerio de Educaci\'on Superior, Ciencia y Tecnolog\'ia (MESCyT) grant FONDOCYT 2022-1D2-091 "Aproximaci\'on Semicl\'asica de Sistemas Hamiltonianos con potenciales Homog\'eneos de dos grados de libertad y osciladores no aut\'onomos".

JTL has been supported by the Spanish State Research Agency (AEI), through the Severo Ochoa and Mar\'ia de Maeztu Program for Centers and Units of Excellence in R\&D (CEX2020-001084-M). Moreover, he has been funded by the Spanish project PGC2018-098676-B-100 funded by MCIN/AEI/10.13039/501100011033  “ERDF A way of making Europe”; by the Spanish project PID2021-122954NB-I00 funded by MCIN/AEI/10.13039/501100011033/ and "ERDF A way of making Europe"; and by a grant from the "Ayudas para la recualificaci\'on del sistema universitario espa\~nol para 2021-2023". He also thanks Laboratorio Subterr\'aneo de Canfranc (LSC) and the Instituto de Biolog\'ia Integrativa de Sistemas (I$^2$Sysbio, CSIC-UV) for their kind hospitality as two of the hosting institutions of this grant.

JJMR has been supported by the Universidad Polit\'ecnica de Madrid research group \emph{Modelos Matem\'aticos no lineales}. As well as he thanks to the Japan Society for Promotion of the Science and to the hospitality of the Departament of Applied Mathematics and Physics of Kyoto University, where he stayed in the fall of 2022 under the Grant  "FY JSPS Invitational Fellowships for Research in Japan" and where part of his contribution to this work was done.

CP is partially supported by the Ministerio de Ciencia e Innovaci\'on grant (PID2019-104658GB-I00); is funded partially by the grant PID-2021-122954NB-100 funded by MCIN/AEI/ 10.13039/501100011033 and
by “ERDF A way of making Europe''.

PA-H and JJMR thanks EPSEB (UPC) for its hospitality during their visit at the final stages of this work.
\newpage
\begin{center}
{\Large{\textsc{Appendix}}    }
\end{center}

\appendix
\section{Bessel, Legendre and Lam\'e equations}
We provide some additional theoretical background to make the paper more self-contained. They aim to complement the Applications section~\ref{appl}.
\subsection{Bessel equation}
\label{apbes}
The well-known Bessel differential equation (see, for instance,~\cite{watson1922}) is a second order linear ode of type:
\begin{equation}\label{beseqldo}
x^2\frac{d^2z}{dx^2}+x\frac{dz}{dx}+\left(x^2-\alpha^2\right) z=0,\qquad \alpha\in\Cset.
\end{equation}
It is well known that for values of the parameter $\alpha\in\mathbb{Z}+\frac{1}{2}$ (i.e. $\alpha$ a half-integer) this equation is integrable in the Picard-Vessiot sense, that is, it admits Liouvillian solutions (see~\cite{morales2013} and references therein). Let us assume that we are in such case, and so
$\alpha\in\mathbb{Z}+\frac{1}{2}$. The general solution of~\eqref{beseqldo} can be expressed as
\[
z(x)=c_1 J_{\alpha}(x)+c_2 Y_{\alpha}(x),
\]
where $J_\alpha$ and $Y_\alpha$ are the (so-called) Bessel functions of first and second kind of parameter $\alpha$, respectively (see~\cite{watson1922}).
The change of variables $y=\sqrt{x}z$ (see~\cite{k1986}) brings the Bessel equation~\eqref{beseqldo} into its normal form, namely,
\begin{equation}\label{redbes1}
y''=r(x)y,\qquad \quad  r(x)=\frac{\alpha^2-\frac{1}{4}}{x^2}-1,
\end{equation}
where $'$ denotes differentiation with respect to $x$. Moreover,
since $\alpha \in\Zset+\frac{1}{2}$ we can set $\alpha=n +\frac{1}{2}$, with $n\in \Zset$. Thus, $r(x)$ above becomes
\[
r(x)=\frac{(\alpha-\frac{1}{2})(\alpha+\frac{1}{2})}{x^2}-1=\frac{n(n+1)}{x^2}-1.
\]
Finally, the transformation $x\mapsto \mu x$, with $\mu \neq 0$, leads equation~\eqref{redbes1} into
\begin{equation}\label{eqrfb}
y''=\left(\frac{n(n+1)}{x^2}-\mu^2\right)y,\qquad n\in\Zset, \qquad \mu\in\Cset^*.
\end{equation}
The case $n=0,-1$ corresponds to a differential equation with constant coefficients, so clearly integrable.
%Furthermore, in the case $\mu=0$ it becomes a Cauchy-Euler differential equation, transformable in a constant coefficients differential equation and, hence, integrable as well.

Equation~\eqref{eqrfb} is very well known in mathematical physics because it corresponds to a radial free particle the Schr\"odinger equation.

\subsection{Legendre equation}\label{A2}
The general Legendre ordinary differential equation takes the form
\begin{equation}\label{LEAp}
(1-z^2)\frac{\partial \xi}{\partial z^2}-2z\frac{\partial \xi}{\partial z}+\left(n(n+1)-\frac{m^2}{1-z^2}\right)\xi=0.
\end{equation}
It has regular singularities at the points $z=-1,1,\infty$. Both singular points $z=\pm 1$ have indices $m/2, -m/2$, while the singular point $z=\infty$ has indices $-n, n+1$.
Hence the differences between the indices at the points $z=-1,1,\infty$ are, up to change in sign, $(m,m,-2n-1)$, respectively. The following Proposition states the values of the parameters $m,n$ for which it is integrable.
It is a correction, completing some missing cases, of~\cite[Proposition 4.]{almp2015}, and it is a direct consequence of Kimura's Theorem~\cite{k1969}.
\begin{proposition}\label{propLegendre}
The Legendre equation~\eqref{LEAp} is integrable if and only if,
either one of the following cases holds:
\begin{itemize}
\item[(A)]  exactly one of the following situations is satisfied
\begin{enumerate}
    \item $n\in\mathbb{Z}$
    \item $m+n\in\mathbb{Z}$, $m\notin \mathbb{Z}$ and $n\notin \mathbb{Z}$
    \item $m-n\in\mathbb{Z}$ $m\notin \mathbb{Z}$ and $n\notin \mathbb{Z}$.
\end{enumerate}

\item[(B)] $m$, $ n$ belong to one of the following families of cases:
\arraycolsep=1.4pt
\def\arraystretch{1.5}
\[
\begin{array}{|c|c|c|c|} \hline
\quad \mbox{Case} \quad & m\in   & n\in  \\ \hline
1 & \quad \frac{1}{2}(2\mathbb{Z} +1) \quad   & \quad \mathbb{C} \quad                 \\[1.2ex]  \hline
%(2) & \mathbb{C}  &   \frac{1}{4}(2\mathbb{Z} +1)              \\[1.2ex]  \hline
2 & \frac{1}{3}(3\mathbb{Z} \pm 1)& \frac{1}{4}(2\mathbb{Z}+1)  \\[1.2ex]  \hline
3 & \frac{1}{3}(3\mathbb{Z} \pm 1)& \frac{1}{6}(6\mathbb{Z}\pm 1)  \\[1.2ex]  \hline
4 & \frac{1}{4}(4\mathbb{Z} \pm 1)& \frac{1}{6}(6\mathbb{Z}\pm 1)  \\[1.2ex]  \hline
5 & \frac{1}{3}(3\mathbb{Z} \pm 1) & \frac{1}{10}(10\mathbb{Z}\pm 3) \\[1.2ex]  \hline
6 & \frac{1}{5}(5\mathbb{Z}\pm 1) & \frac{1}{6}(6\mathbb{Z}\pm 1)  \\[1.2ex]  \hline
7 & \frac{1}{5}(5\mathbb{Z}\pm 2) & \frac{1}{10}(10\mathbb{Z}\pm 3)  \\[1.2ex] \hline
8 & \frac{1}{5}(5\mathbb{Z}\pm 1) & \frac{1}{10}(10\mathbb{Z}\pm 1)  \\[1.2ex] \hline
\end{array}
\]
\end{itemize}
\end{proposition}
%In this appendix we present the complete proof of Proposition \ref{propLegendre} and is a correction of Proposition 4 in~\cite{almp2015}.
%Legendre's equation~\eqref{LE3} has the regular singularities $\pm 1$  both with the same exponents  $\pm m/2$ and the regular singularity $\infty$ has exponents $-n$ and $n+1$.
Let us deal with the proof of this proposition. We denote by $\lambda=m$, $\mu=m$ and $\nu=-(2n+1)$ the exponent differences
and we analyze all the cases of Kimura's Theorem \cite{k1969}. Hence, in order the Legendre equation~\eqref{LE3} to have Liouvillian solutions it is necessary and sufficient that, either case A or case B holds. Let us proceed case by case.

{\bf Case A}. At least one of  $\lambda+\mu+\nu$, $-\lambda+\mu+\nu$, $\lambda-\mu+\nu$ or $\lambda+\mu-\nu$  is an odd integer. Equivalently, at least one of the following relations is satisfied:
$\lambda+\mu+\nu\in 2\mathbb{Z}+1$, $-\lambda+\mu+\nu\in 2\mathbb{Z}+1$, $\lambda-\mu+\nu\in 2\mathbb{Z}+1$ or $\lambda+\mu-\nu\in 2\mathbb{Z}+1.$

We consider each item separately.
\begin{description}
\item[1.1] Relation $\lambda+\mu+\nu=2m-2n-1\in 2\mathbb{Z}+1$ yields to $m-n\in\mathbb{Z}$. In conclusion $(m,n)\in \mathbb{Z}^2$ or $m\notin \mathbb{Z}$ and $n\notin \mathbb{Z}$ with $m-n\in \mathbb{Z}$.
\item[1.2] Relation $-\lambda+\mu+\nu=2n-1\in 2\mathbb{Z}+1$ yields to $n\in\mathbb{Z}$. In conclusion $n\in \mathbb{Z}$ and $m\in\mathbb{C}$.
\item[1.3] Relation $\lambda-\mu+\nu= 2n-1\in 2\mathbb{Z}+1$ yields to $n\in\mathbb{Z}$. In conclusion $n\in \mathbb{Z}$ and $m\in\mathbb{C}$.
\item[1.4] Relation $\lambda+\mu-\nu=2m+2n+1\in 2\mathbb{Z}+1$ yields to $m+n\in\mathbb{Z}$. In conclusion $(m,n)\in \mathbb{Z}^2$ or $m\notin \mathbb{Z}$ and $n\notin \mathbb{Z}$ with $m+n\in \mathbb{Z}$.
\item[1.5] Any combination of the previous four relations does not provide any new condition on the parameters $n$ and $m$.
\end{description}
In this way, the first part of Kimura's Theorem for Legendre differential equation~\eqref{LE3} is proved and, summarising, we obtain:
\begin{enumerate}
\item $n\in\Zset$; %$m\in \Zset$
%\item $n\in\mathbb{Z}$, $m\notin \Zset$
\item $m+n\in\Zset$, $m\notin \Zset$ and $n\notin \Zset$;
\item $m-n\in\Zset$ $m\notin \Zset$ and $n\notin \Zset$;
\end{enumerate}

\medskip

{\bf Case B}. The quantities $\lambda$ or $-\lambda$, $\mu$ or $-\mu$, $\nu$ or $-\nu$ take, in an arbitrary order, values given in Kimura's table, see \cite{k1969}.
The cases 4, 6, 9, 10, 12, 14 and 15 of Kimura's table are discarded because for Legendre equation~\eqref{LE3} we have that two differences of exponents have the value.

Now we check the rest of the cases of Kimura's table:
\begin{description}

\item[Case 1.]  By case 1 in Kimura's table, we have $\lambda=\mu=m\in\pm (\Zset+\frac{1}{2})= (\Zset\pm\frac{1}{2})=\frac{1}{2}(2\Zset\pm 1)$ and $\nu=-(2n+1)\in\mathbb{C}$, which lead us to $m\in\frac{1}{2}(2\Zset+ 1)$ and $n\in\mathbb{C}$.
Thus, we obtain the case 1 in the table provided in Proposition \ref{propLegendre}.

\item[Case~2.]  By case 2 in Kimura's table %and due to $\lambda=\mu=m$,
the only one possibility for $m$ and $n$ is provided by $\lambda=\mu=m\in\pm (\Zset+\frac{1}{3})= (\Zset\pm\frac{1}{3})=\frac{1}{3}(3\Zset\pm 1)$ and $\nu=-(2n+1)\in\pm (\Zset+\frac{1}{2})= (\Zset\pm\frac{1}{2})=\frac{1}{2}(2\Zset\pm 1)$, which lead us to $m\in\frac{1}{3}(3\Zset\pm 1)$ and $n\in \frac{1}{4}(2\Zset+1)$. Thus, we obtain the case 2 in the table provided in Proposition \eqref{propLegendre}.

\item[Case~3.]  By case 3 in Kimura's table % and due to the even condition
%\tlcom{Cual es la 'even condition'?}}
%$\lambda=\mu=m$,
%$q\in 2\Zset$,
the only possibility for $m$ and $n$ is provided by $\lambda=\mu=m\in\pm (\Zset+\frac{1}{3})= (\Zset\pm\frac{1}{3})=\frac{1}{3}(3\Zset\pm 1)$ and $\nu=-(2n+1)\in\pm (2\Zset+\frac{2}{3})= (2\Zset\pm\frac{2}{3})$. Therefore $n+\frac{1}{2}\in \Zset\pm \frac{1}{3}$, which leads to $m\in\frac{1}{3}(3\Zset\pm 1)$ and $n\in \frac{1}{6}(6\Zset\pm 1)$. Thus, we obtain the
case~3 at the table provided in the statement of the Proposition.

\item[Case~5.]  By case~5 in Kimura's table and due to the even condition
%$\lambda=\mu=m$,
%$q\in 2\Zset$
the only possibility for $m$ and $n$ is provided by $\lambda=\mu=m\in\pm (\Zset+\frac{1}{4})= (\Zset\pm\frac{1}{4})=\frac{1}{4}(4\Zset\pm 1)$ and $\nu=-(2n+1)\in\pm (2\Zset+\frac{2}{3})= (2\Zset\pm\frac{2}{3})$. Therefore $n+\frac{1}{2}\in \Zset\pm \frac{1}{3}$, which lead us to $m\in\frac{1}{4}(4\Zset\pm 1)$ and $n\in \frac{1}{6}(6\Zset\pm 1)$. Thus, we obtain the
case~4 in the table of the statement.

\item[Case~7.]  By case~7 in Kimura's table and due to the even condition %$\lambda=\mu=m$
  %$q\in 2\Zset$,
the only possibility for $m$ and $n$ is provided by $\lambda=\mu=m\in\pm (\Zset+\frac{1}{3})= (\Zset\pm\frac{1}{3})=\frac{1}{3}(3\Zset\pm 1)$ and $\nu=-(2n+1)\in\pm (2\Zset+\frac{2}{5})= (2\Zset\pm\frac{2}{5})$. Therefore $n+\frac{1}{2}\in \Zset\pm \frac{1}{5}$, which lead to $m\in\frac{1}{3}(3\Zset\pm 1)$ and $n\in \frac{1}{10}(10\Zset\pm 3)$. Thus, we obtain the case~5 in the table provided in Proposition \ref{propLegendre}.

\item[Case~8.]  By case~8 in Kimura's table and due to the even condition %$\lambda=\mu=m$,
%$q\in 2\Zset$
the only possibility for $m$ and $n$ is provided by $\lambda=\mu=m\in\pm (\Zset+\frac{1}{5})= (\Zset\pm\frac{1}{5})=\frac{1}{5}(5\Zset\pm 1)$ and $\nu=-(2n+1)\in\pm (2\Zset+\frac{2}{3})= (2\Zset\pm\frac{2}{3})$. Therefore $n+\frac{1}{2}\in \Zset\pm \frac{1}{3}$, which leads to $m\in\frac{1}{5}(5\Zset\pm 1)$ and $n\in \frac{1}{6}(6\Zset\pm 1)$. Thus, we obtain the
case~6 in the table provided in the Proposition.

\item[Case~11.]  By case 11 in Kimura's table and due to the even condition %$\lambda=\mu=m$,
 %$q\in 2\Zset$
the only possibility for $m$ and $n$ is provided by $\lambda=\mu=m\in\pm (\Zset+\frac{2}{5})= (\Zset\pm\frac{2}{5})=\frac{1}{5}(5\Zset\pm 2)$ and $\nu=-(2n+1)\in\pm (2\Zset+\frac{2}{5})= (2\Zset\pm\frac{2}{5})$. Therefore $n+\frac{1}{2}\in \Zset\pm \frac{1}{5}$, which leads to $m\in\frac{1}{5}(5\Zset\pm 2)$ and $n\in \frac{1}{10}(10\Zset\pm 3)$. Thus, we obtain the case~7 in the table.

\item[Case~13.]  By case~13 in Kimura's table and due to the even condition %$\lambda=\mu=m$,
  %$q\in 2\Zset$ t
the only possibility for $m$ and $n$ is provided by $\lambda=\mu=m\in\pm (\Zset+\frac{1}{5})= (\Zset\pm\frac{1}{5})=\frac{1}{5}(5\Zset\pm 1)$ and $\nu=-(2n+1)\in\pm (2\Zset+\frac{4}{5})= (2\Zset\pm\frac{4}{5})$. Therefore $n+\frac{1}{2}\in \Zset\pm \frac{2}{5}$, leading to $m\in\frac{1}{5}(5\Zset\pm 1)$ and $n\in \frac{1}{10}(10\Zset\pm 1)$. Thus, we obtain the case 8 in the table.   \end{description}

\subsection{Lam\'e Equation}\label{apLame}
The well-known Lam\'e ordinary differential equation is, in Weierstrass form, given by
\begin{equation}\label{Lame}
\frac{d^2y}{dz^2}=(n(n+1)\wp(z)+B)y,
\end{equation}
where $\wp$ is the Weierstrass function, a solution of the differential equation $(dw/dz)^2=h(w)$, with $h(w)=4w^3-g_2w-g_3$ and discriminant $\Delta=g_2^3-27g_3^2\neq 0$ (in order the polynomial $h(w)$ to have simple roots; otherwise it could be transformed into a simpler form).
%The $\zeta$ function satisfies $\zeta'=\wp$, i.e,, $\zeta$ Weierstrass function is the antiderivative of $\wp$ Weierstrass function.
This equation~\eqref{Lame} depends on four parameters: $n$, $B$, $g_2$ and $g_3$.

We assume the basic periods of $\wp$, named $2\omega_1$, $2\omega_3$, to be real and purely imaginary, respectively.
These conditions are satisfied when $g_2$ and $g_3$ are real and $\Delta>0$. If we denote by $e_1$, $e_2$ and $e_3$ the roots of $h(w)$, we have that they are all real and it is no restrictive to assume that $e_3<e_2<e_1.$
Then $\omega_i=e_i$, $i=1,2,3$, being $\omega_1$ real and $\omega_3$ purely imaginary. The function $g(x)=\wp(x+\omega_3)$ is real and regular for $x\in \Rset$ (see for instance~\cite{ww2020}).

Some integrability conditions of Lam\'e equation \eqref{Lame}, in the sense of the differential Galois theory, can be found in~\cite{ms1996, morales2013, ww2020}. They are:
 %In particular for $n\in\mathbb{N}$ and $B$, arbitrary, $g_2$ and $g_3$,  Lam\'e
%equation becomes integrable in three cases, see \cite{morales2013}.

\begin{itemize}
\item[(i)] {\it{Lam\'e and  Hermite-Halphen solutions}}, $n\in\Nset$.

\item[(ii)] {\it{Brioschi--Halphen--Crawford solutions}}, $n+1/2 \in \Nset$ and some algebraic conditions on the rest of the parameters.

\item[iii)] {\it{Baldasarri solutions}}, $n+\frac{1}{2}\in\frac{1}{3}\Zset\bigcup \frac{1}{4}\Zset \bigcup \frac{1}{5}\Zset\setminus\Zset$ and some other involved conditions on the rest of parameters.
\end{itemize}
In particular, for $n\in \Nset$, Lam\'e equation becomes integrable.

%\printbibliography
%\bibliography{ALMP}{}
%\bibliographystyle{plain}
%\bibliographystyle{unsrt}
%\bibliography{ALMP}

\end{document}